\documentclass[11pt]{amsart}

\usepackage{amsmath,amsthm,amsfonts,amssymb,url}
\usepackage[T1]{fontenc}
\usepackage{xypic}

\newcommand{\Q}{\mathbb{Q}}
\newcommand{\Z}{\mathbb{Z}}

\newcommand{\Frob}{{\rm{Frob}}}

\newcommand{\U}{\mathrm{U}}

\newcommand{\Gal}{\mathrm{Gal}}

\newtheorem{theorem}{Theorem}

\newtheorem{prop}{Proposition}

\newtheorem{lemma}{Lemma}

\theoremstyle{definition}
\newtheorem{remark}{Remark}

\def\]{\textup{\mbox{]\hspace{-.15em}]}}}
\def\[{\textup{\mbox{[\hspace{-.15em}[}}}

\newcommand{\anneau}{\mathcal{O}}

\newcommand{\Li}{\text{Li}\,}
\date{}

\title{Remarks on the error term in Chebotarev's density theorem}

\begin{document}
\baselineskip 17pt

\address{Jo\"el Bella\"iche\\Brandeis University\\
415 South Street\\Waltham, MA 02454-9110\\U.S.A}
\email{jbellaic@brandeis.edu}
\author{Jo\"el Bella\"iche}
\maketitle

\par \bigskip
\par \bigskip

\section{The error term in Chebotarev's density theorem}

In \cite{effective}, and then again, without the third author, in \cite{eff}, Ram Murty, Kumar Murty and Saradha prove a result (conditional to some standard conjectures of analytic number theory, namely the General Riemann Hypothesis and the Artin conjecture) bounding the error term in Chebotarev's Density Theorem (see (\ref{mms}) below) and ask whether an improved bound (see (\ref{mmsq}) below) holds. The aim of this article is to prove that the answer to this question is negative, to answer (by the affirmative) a related question of Serre from \cite{serre}, and more generally, to discuss the best possible error terms (of the forms considered in \cite{effective} and \cite{eff}) in Chebotarev's density theorem. 

\par \bigskip

To formalize this question of Ram Murty, Kumar Murty and Saradha precisely, let us fix some notations:
in all this paper, $L$ will denote a finite Galois extension of $\Q$ of degree $n$ and Galois group $G$, and $D$ will denote a conjugacy set (that is, a union of conjugacy 
classes) in $G$. We will denote by $M$ the product of all prime numbers that are ramified in $L$. For $x$ a positive real number, we call $\pi_D(x)$ the number of primes $p < x$ such that the conjugacy class $\Frob_{p,L/\Q}$ of $G$ is contained in $D$. 
By (GRH) we shall mean, as in \cite{effective} the Generalized Riemann Hypothesis for all Artin $L$-functions.

The first effective version of Chebotarev's density theorem, proved by Lagarias and Odlyzko (\cite{lo}) 
and soon thereafter improved by Serre (\cite[($20_R$), page 134]{serre}) states, that, assuming (GRH):
\begin{eqnarray} \label{los1} \pi_D(x) = \frac{|D|}{|G|} \Li(x)+O\left( x^{1/2} |D| (\log x + \log |G| + \log M)\right). \end{eqnarray}
The implied constant in the above formula is absolute. In this paper we shall not be interested in the logarithmic terms in $x$ and $|G|$. Therefore, let us replace
(\ref{los1}) by its following slight weakening: 
\begin{eqnarray} \label{los2} \text{For every $\epsilon>0$, \ \ \ } \pi_D(x) = \frac{|D|}{|G|} \Li(x) + O\left(x^{1/2+\epsilon} |D| |G|^{\epsilon} \log M \right). \end{eqnarray}
In this formula the implied constant depends only on $\epsilon$.

Later, Murty, Murty and Saradha proved (cf. \cite[Corollary 3.7]{effective}), using (GRH) and the Artin's conjecture of holomorphy of Artin's $L$-functions, that 
\begin{eqnarray} \label{mms} \pi_D(x) = \frac{|D|}{|G|} \Li(x) + O\left(x^{1/2} |D|^{1/2} (\log x + \log |G| + \log M) \right) \end{eqnarray}
with an absolute implied constant, which again can be slightly weakened into:
\begin{eqnarray} \label{mms2} \text{For every $\epsilon>0$, \ \ \ } \pi_D(x) = \frac{|D|}{|G|} \Li(x) + O\left(x^{1/2+\epsilon} |D|^{1/2} |G|^{\epsilon} \log M \right) 
\end{eqnarray}
%

In the same paper \cite[\S3.13]{effective} (and also in \cite{eff} without the third author), Murty, Murty and Saradha ask the following question. Let $\alpha(G)$ be the number of conjugacy classes in $G$.
Is it true that (with an absolute implied constant)
\begin{eqnarray} \label{mmsq} \pi_D(x) = \frac{|D|}{|G|} \Li(x) + O\left(x^{1/2+\epsilon} |D|^{1/2} \alpha(G)^{-1/2} (\log x + \log |G| + \log M) \right)\ ? \end{eqnarray}
 Once again, a positive answer to this question would mean a positive answer to its following simplified version. Is it true that (with an implied constant depending only on $\epsilon$)
\begin{eqnarray} \label{mmsq2} \text{For every $\epsilon>0$, \ \ \ } \pi_D(x) = \frac{|D|}{|G|} \Li(x) + O\left( x^{1/2+\epsilon} |D|^{1/2} \alpha(G)^{-1/2} 
|G|^{\epsilon} \log M\right) \ ?
\end{eqnarray} 
Here the implied constant depends only on $\epsilon$.

The answer to (\ref{mmsq}) and (\ref{mmsq2}) is no for a trivial reason: {\it there is no specified range}. Indeed, assume that $|D|=1$ to fix ideas, and that $M$ is constant. 
If $x < \alpha(G)^c$ for some constant $c<1$, then the error term in those formulas goes to $0$ when $\alpha(G)$ goes to infinity, 
hence $\pi_D(x+1)-\pi_D(x)$ goes to $0$ when $x$ and $\alpha(G)$ go to infinity with $x<\alpha(G)^c$. So, for $\alpha(G)$ large enough and
$x < \alpha(G)^c$,  $\pi_D(x+1)=\pi_D(x)$ since those numbers are integers.
Now, choose a field $L/\Q$ unramified outside $M$, with Galois group $G$ such that $\alpha(G)$ is large enough in the preceding sense (of course there are plenty of such field even among cyclotomic fields) and take $D=\{\Frob_p\}$ for a prime $p < \alpha(G)^c$. Then clearly $\pi_D(p+1) = \pi_D(p)+1$, a contradiction.

We note that this omission of the range made in  \cite{effective} and in \cite{eff} seems quite frequent in similar questions 
in the literature of analytic number theory, for example those concerning the primes in an arithmetic progression with common difference $q$: cf.  \cite[First paragraph of page 309]{m2} or \cite[(17.5)]{ik}, which are both false as stated  for the same trivial reason.
However, the remaining of the discussion in \cite[\S17.1]{ik} makes clear that the restriction $x>q$ is implicitly assumed, and it is likely that such a restriction was
also in the authors' mind in \cite{m2}.

Similarly, it is natural in the questions (\ref{mmsq}) and (\ref{mmsq2}) to restrict our attention 
to the range  $ x \geq |G| $ or even, to be more prudent, to $x \geq |G| \log^\alpha |G|$
for every $\alpha>0$. That is to say, one is led to ask:
\begin{eqnarray} \label{mmsq3} & \text{  For every $\epsilon>0$, is it true that if $x>|G|$,} \nonumber \\ & \pi_D(x) = \frac{|D|}{|G|} \Li(x) + O\left( x^{1/2+\epsilon} |D|^{1/2} \alpha(G)^{-1/2} 
|G|^{\epsilon} \log M\right) ,\\ \nonumber
& \text{the implied constant depending only on $\epsilon$}?
\end{eqnarray} 
 And:
\begin{eqnarray} \label{mmsq4}\nonumber & \text{For every $\epsilon>0$ and every $\alpha>0$, is it true that if $x>|G| (\log |G|)^\alpha$,} \\ & \pi_D(x) = \frac{|D|}{|G|} \Li(x) + O\left( x^{1/2+\epsilon} |D|^{1/2} \alpha(G)^{-1/2} 
|G|^{\epsilon} \log M\right),\\
\nonumber & \text{the implied constant depending only on $\epsilon$ and $\alpha$}?
\end{eqnarray}
\begin{remark}
Consider the cyclotomic case $L=\Q(\mu_q)$ (so $G=(\Z/q\Z)^\ast$) and $D=\{d\}$, where $q$ and $d$ are relatively prime positive integers. In this case, $\pi_D(x)$ counts the primes $p<x$ such that $p \equiv d \pmod{q}$, and a conjecture of
Friedlander and Granville (\cite[Conjecture 1(b), page 366]{fg}), correcting the famous conjecture of Montgomery (\cite{m1}, \cite{m2}), states
\begin{eqnarray} \nonumber & \text{For every $\epsilon>0$, if $x>q$ }  \\  \label{fgconj} & \pi_D(x) = \frac{1}{\phi(q)} \Li(x) + O\left( x^{1/2+\epsilon} q^{-1/2}
\right), \end{eqnarray} 

Using the well-known estimate $q/\log \log q < \phi(q) < q$, it is an easy exercise to see that, in this cyclotomic case with $|D|=1$,
$ (\ref{fgconj}) \Longrightarrow (\ref{mmsq4}) $. At any rate it is plain that
 (\ref{mmsq3}) and (\ref{mmsq4}) {\it in the cyclotomic case with $|D|=1$} on the one hand, and (\ref{fgconj}) on the other hand are very close, and 
 it is almost certain that those three conjectures hold or fall together.
\end{remark}

\par \bigskip
Back to general case, our main result is that, even if we add the forgotten restriction of the range, the answer to Murty, Murty, and Seradha's question is no.
\begin{theorem} The answer to (\ref{mmsq4}), hence to (\ref{mmsq3}), is no.
\end{theorem}
We will prove this theorem twice, by giving two separate counter-examples: one in the case $G$ dihedral, $D=\{1\}$ (see Propostion~\ref{dihedral}), the second in the case $|G|$ abelian, $|D|$ large (see Proposition~\ref{abelian}). The fact that two natural generalizations of the conjecture of Friedlander and Granville fail suggest that this conjecture is either extremely difficult (wihich of course, was clear to begin with), or even, possibly, false.

 \par \bigskip
 Let us now discuss the question of Serre. A result of Lagarias and Odlyzko (\cite{lo2}) states that, under (GRH),  $\pi_D(x) > 0$ if $x> c (\log d_K)^2$, where $c$
 is some absolute constant. Serre asks (\cite[Remark 2, page 135]{serre}) whether the exponent $2$ in this formula is the best possible.
 In view of \cite[Proposition 6]{Serre}, which in our notations states $$|G| \log M / 2 \leq \log d_K \leq (|G|-1) \log M + |G| \log |G|,$$ the question is equivalent to: 
\begin{eqnarray} \label{jpsq}\ \ \ \ \text{ Does there exist $e<2$ and $c>0$, such that $x> c |G|^e \log M$ implies $\pi(x)>0$?} \end{eqnarray}
We shall see (cf. Prop.~\ref{dihedral}) that the answer to (\ref{jpsq}) is no -- that is, the answer to Serre's question "is $2$ the best possible constant?" is yes.

\par \bigskip 

On a more constructive tone, one may then ask: what is the best possible error term in Chebotarev's density theorem? 
Of course, in this generality, the question is not interesting (the answer is obviously $|\pi_D(x) - \frac{|D|}{|G|} \Li(x)|$): we need to restrict our study to
some special families of allowed error terms, for example error terms depending only on the size of the group $G$ and the set $D$ through the product of a power of $|G|$ and a power of $|D|$.

\begin{theorem}
For two real numbers $a$ and $b$, consider the following assertion:
 
\begin{eqnarray*} &  \text{For every $\epsilon, \alpha>0$, in the range $x > |G| (\log |G|)^\alpha$,} \\  (C_{a,b})\ \ \ \ \ \  &  \pi_D(x) = \frac{|D|}{|G|} \Li(x) + O\left(x^{1/2+\epsilon} |D|^a |G|^{b+\epsilon} \log M \right)  
\\ & \text{with an implied constant depending only on $a$, $b$, $\epsilon$ and $\alpha$.}
\end{eqnarray*}

Then for $(C_{a,b})$ to be true it is necessary, and, under (GRH), sufficient that $b \geq 0$ and $a+b \geq 1/2$.
\end{theorem}
In other words, the estimate of Ram Murty, Kumar Murty and Saradha (\ref{mms2}), or $(C_{1/2,0})$ in our language, 
is the best possible among estimates of the form $(C_{a,b})$. Indeed,  
if $a+b \geq 1/2$,  and $b \geq 0$, $|D|^{1/2} \leq |D|^{a+b} = |D|^a |D|^b  \leq |D|^a |G|^b$ because $|D|<|G|$, and $(C_{1/2,0})$ implies $(C_{a,b})$,
hence the necessity in the above theorem. 
The sufficiency will be proved in Prop.~\ref{dihedral} and Prop.~\ref{abelian}.

One may then ask for other types of error estimates. For example, question (\ref{mmsq2}) suggests that we look at error term of the form $O\left(x^{1/2+\epsilon} |D|^a \alpha(G)^b |G|^{\epsilon} \log M \right)$. 
\begin{theorem}
For two real numbers $a$ and $b$, consider the following assertion:
 
\begin{eqnarray*} 
 &  \text{For every $\epsilon>0$, $\alpha>0$,   in the range $x > |G| (\log |G|)^\alpha$,} \\  (C'_{a,b})\ \ \ \ \ \   & \pi_D(x) = \frac{|D|}{|G|} \Li(x) + O\left(x^{1/2+\epsilon} |D|^a \alpha(G)^{b} |G|^\epsilon \log M \right) \\
& \text{with an implied constant depending only on $a$, $b$, $\epsilon$ and $\alpha$.}\end{eqnarray*}
Then for $(C'_{a,b})$ to be true it is necessary that $b \geq 0$ and $a+b \geq 1/2$.
\end{theorem}
In this case, since $\alpha(G)$ may be larger or smaller than $|D|$, it is not clear in this case that $(C'_{1/2,0})$ (which is the same as $(C_{1/2,0})$, hence a theorem under (GRH)) implies $(C'_{a,b})$ if $b \geq 0$, $a+b \geq 1/2$. That is to say, we ignore (even under (GRH))  if $(C'_{a,b})$ is true  for some $b>0$, $a+b=1/2$, for example for $a=1/4$, $b=1/4$. This would be very surprising however, as it would mean a better bound, for the same $|G|$ and $|D|$, in the non-abelian case than in the abelian case, and at any rate contrary to the intuition underlying Question 
(\ref{mmsq2}) where $\alpha(G)$ appeared with a negative exponent $b$, hence giving a better bound in the abelian than in the non-abelian case. 

Let us end this already too long introduction that we can get better and useful estimate for the error term which make intervene the {\it structure} of the group $G$ and the set $D$ rather than just their size $|G|$ and $|D|$. See \cite[page 143]{ik} and \cite{B}.

\section{A family of dihedral examples with $D=\{1\}$}

The examples will depend on an integral parameter $n=2^r$ with $r \geq 2$.
Let $\anneau$ be the order $\Z[\sqrt{ - n^2}] = \Z[n i]$ in the field of Gaussian numbers $\Q(i)$. 
Let $L$ be the ring class field of $\anneau$. As usual let $G=\Gal(L/\Q)$, $M$ the product of all primes ramified in $L$. We set $D=\{1\}$, so $\pi_D(x)$ counts the primes less than $x$ that are totally split in $L$.

\begin{prop} \label{dihedral} 
For $L$, $D$ as above and $n = |G|$ large enough, $(C_{a,b})$ and $(C'_{a,b})$ are false whenever $b<0$. Moreover (\ref{jpsq}) is false whenever $e<2$.
\end{prop} 

We collect in the following lemma various elementary results, mainly from Cox's book.
\begin{lemma} \label{lemmad}
\item[(i)] $M=2$
\item[(ii)] $G$ is a dihedral group of order $n$.
\item[(iii)] If $x = |G|^2=n^2$, then $\pi_D(x)=0$.
\item[(iv)] $\alpha(G)>n/4$. 
\end{lemma}
\begin{pf}  The conductor $f = [\Z[i] : \anneau]$ of the order $\anneau$ is $n$, hence the class number of that order is $h(\anneau)=\frac{f}{[\Z[i]^\ast:\anneau^\ast]}=\frac{f}{2}  = n/2 =2^{r-1}$ according to \cite[Theorem 7.24]{cox}. Hence by definition of the ring class field and by class field theory, $L$ is an abelian extension of $\Q(i)$ of degree $n/2$, ramified only at (the prime dividing) $2$. Therefore, $G$ is dihedral of order $n$, and $M=2$. This proves (i) and (ii).

By \cite[Theorem 9.4]{cox}, an odd prime $p$ is totally split in $L$ if and only if it is represented by the form  
$$a^2 + n^2 b^2.$$ 
It is clear that no prime $p \leq n^2$ is of this form, hence $\pi_D(x)=0$ if $x<n^2$, which is (iii).

Finally (iv) follows from the computation of the number of conjugacy classes in a dihedral group, which is easy and standard (see e.g. \cite[\S5.3]{serrerep})
\end{pf}
\par \bigskip

Now assume that formula ($C_{a,b}$) or ($C'_{a,b})$ is true in the case under consideration for some value of the parameters $a$ and $b$ such that $b<0$.
We are going to obtain a contradiction.  Since $\alpha(G) < |G|$, it suffices to consider the case of $C'_{a,b}$:
$$\pi_D(x) = \frac{1}{n} \Li(x) + O(x^{1/2+\epsilon} \alpha(G)^b n^\epsilon).$$
By Lemma~\ref{lemmad}(iv), $\alpha(G)^b = O (n^b)$, hence 
 $$\pi_D(x) = \frac{1}{n} \Li(x) + O(x^{1/2+\epsilon} n^{b+\epsilon}).$$
 Let us fix some $\alpha>0$. If $n$ is large enough, $x=n^2$ is certainly in the range $x > n \log(n)^\alpha$, so $\pi_D(x)=0$, and the formula becomes
 $\frac{1}{n} \Li(n^2) = O(n^{1+b+3 \epsilon})$ or $\Li(n^2)=O(n^{2+b+3 \epsilon})$.
 Since $b<0$, we can choose $\epsilon > 0$ such that $-\epsilon' := b+3 \epsilon <0$, and we get $\Li(n^2)=O(n^{2-\epsilon'})$, which is absurd
 since $\Li(n^2) \sim n^2 / (2\log n)$. This completes the proof of the first part of Proposition~\ref{dihedral}.
 
Concerning Serre's question (\ref{jpsq}), one has just seen that $\pi_D(x) =0$ for all $x \leq |G|^2 = \frac{1}{\log 2} |G|^2 \log M$, which shows
that (\ref{jpsq}) cannot be true for $e<2$. This completes the proof.
  
\section{A family of abelian examples with $|D|$ large}
 
Again, the example will depend on a integral parameter $n=2^r$ with $r \geq 2$. We define $L = \Q(\mu_{2^{r+1}})=\Q(\mu(2n))$, so that 
$G=\Gal(L/\Q)=(\Z/2n\Z)^\ast$ has order $n$, and for an odd prime $p$, $\Frob_{p, L/\Q}$ is just $p \pmod{2n}$.  One also have $M=2$, and $\alpha(G)=|G|=n$.
Fix an $\alpha>0$, such that $\alpha<1$. We let $D \subset (\Z/2n\Z)^\ast$  be the set of all odd residue classes $d$
modulo $2n$ such that the arithmetic progression $d + 2n \Z$ does not contain any prime smaller than $n \log(n)^\alpha$. 
Clearly, the complement of $D$ has size at most $\pi(n \log(n)^\alpha) = o(n)$ by the prime number theorem, so $|D| \sim n$ when $n$ goes to infinity.
On the other hand, by definition $\pi_D(n \log(n)^\alpha) = 0$. 

\begin{prop} \label{abelian}  For $L$, $D$ as above and $n = |G|$ large enough, $C_{a,b}$ and $C'_{a,b}$ is false if $a+b < 1/2$
\end{prop}
Indeed, $C_{a,b}$ and $C'_{a,b}$ are identical in this case, and state that for $x \geq n \log (n)^\alpha$,
$$\pi_D(x) = \frac{|D|}{n} \Li(x) + O\left(x^{1/2+\epsilon} n^{a+b+\epsilon} \right) .$$
Applying this to $x = n \log(n)^\alpha$, we get, since $\pi_D(x)=0$,
$$ \frac{|D|}{n} \Li(x) = O (n^{a+b+1/2+ 2 \epsilon} \log(n)^{\alpha(1/2+\epsilon)} ) $$
which implies since $|D| \sim n$, and $n \leq x$, $$n = O (n^{a+b+1/2+2 \epsilon} \log(n)^\beta)$$ where $\beta$ is some real number depending only on $\epsilon$ and $\alpha$. If $a+b < 1/2$, on can choose $\epsilon > 0$ such that the exponent of $n$ in the RHS is $<1$, giving a contradiction.

\end{document}